\documentclass[12pt]{article}
\usepackage{amssymb,amsmath,theorem}
\newtheorem{thm}{Theorem}

{\theorembodyfont{\rmfamily}
\newtheorem{rem}{Remark}}
\def\P{{\mathbb P}}
\def\E{{\mathbb E}}

\def\bbr{{\mathbb R}}

\def\nubar{\overline{\nu}}

\begin{document}
\title{Median, Concentration and Fluctuations  
for L\'evy Processes}
\author{Christian Houdr\'e\thanks{School of Mathematics, Georgia Institute of Technology, 
Atlanta, GA 30332, USA, houdre@math.gatech.edu.  This research was done in part 
while I visited  L'\'Ecole Normale Sup\'erieure 
whose hospitality and support is gratefully acknowledged.}  and Philippe Marchal\thanks{DMA, \'Ecole Normale Sup\'erieure, 
75005 Paris, France, Philippe.Marchal@ens.fr}}
\maketitle

\begin{abstract}
We estimate a median of $f(X_t)$ where $f$ is a Lipschitz function,
$X$ is a L\'evy process and $t$ is an arbitrary time.  This leads to 
concentration inequalities for $f(X_t)$.  In turn, 
corresponding fluctuation estimates are obtained under assumptions 
typically satisfied if the process has a regular behavior in small time 
and a, possibly different, regular behavior in large time.   
\end{abstract}

\section{Introduction}

In $\bbr^d$, let $X = (X_t,t\geq 0)$ be a L\'evy process, without Gaussian component.   
Its characteristic exponent $\psi_X$ is given, 
for all $u\in\bbr^d$, by 
$$
\E\exp(i\langle u,X_t\rangle)=\exp(t\psi_X(u)), 
$$
where 
\begin{equation}\label{eq1}
\psi_X(u)=i\langle u, b \rangle + \int_{\bbr^d} \left(e^{i\langle u,y\rangle}
-1-i\langle u,y\rangle {\bf 1}_{\|y\|\le 1}\right)\nu (dy),
\end{equation}
$b\in\bbr^d$ and $\nu\not\equiv 0$ (the L\'evy measure)
is a positive Borel measure without atom at the origin and such that
$\int_{\bbr^d} (\|y\|^2\wedge 1)\nu (dy)<+\infty$ (throughout,
$\langle\cdot,\cdot\rangle$ and $\|\cdot\|$ are respectively the 
Euclidean inner product and norm in $\bbr^d$).\\

While the asymptotic behavior of $X$ in small or large time can
be deduced from the asymptotic behavior of $\psi$ near the origin or at infinity, 
it is more difficult to get precise estimates, for the law of $X_t$, at some given time 
$t$.  
However, when $X$ has finite mean, Marcus and Rosi\'nski \cite{MR} (see the next section 
for a precise statement) provide a fine estimation of 
$\E\|X_t\|$ involving the functions
$$
V(R)=\int_{\|x\|\leq R} \|x\|^2 \nu(dx), 
$$
and
$$
M(R)=\int_{\|x\|> R} \|x\|  \nu(dx),   
$$ 
$R>0$.

If one removes the assumption of finite mean, in which case $M(R)$ becomes infinite, 
the natural way to express the order of magnitude of $\|X_t\|$ is to consider 
one of its medians.  One may then want to estimate this median and to further know 
how $\|X_t\|$ is concentrated around it.  More generally, one may 
ask the same question for $f(X_t)$, where $f$ is a Lipschitz 
function with respect to the Euclidean norm.  
The aim of this paper is to investigate these questions and related ones.

In essence, the main result of the present paper is that under some rather general 
hypotheses, if $f$ is a Lipschitz function with Lipschitz constant 1 (a 1-Lipschitz function), 
the order of magnitude of the median and of the fluctuations of $f(X_t)$ is
given by functions of the form
\begin{equation}\label{h}
h_c(t)=\inf\left\{x>0: \frac{V(x)}{x^2}=\frac{c}{t}\right\},  
\end{equation}
where $c$ is some positive real. More precisely, 
denote by $\nubar$ the tail of $\nu$, i.e., let  
$$
\nubar(R)=\int_{\|x\|> R} \nu(dx), 
$$
then we have:

\begin{thm}
Let $X$ be a L\'evy process with characteristic exponent \eqref{eq1}.
Let $f$ be a 1-Lipschitz function, let $c>0$, let $t>0$, and 
let $h_c$ be given by \eqref{h}.  Then, for every $t > 0$ such that 
\begin{equation}\label{tnu}
t\nubar(h_c(t)) \leq 1/4,
\end{equation}
any median $mf(X_t)$ of $f(X_t)$ satisfies:
\begin{equation}\label{eqmed}
|mf(X_t)-f(0)|\leq h_c(t)\left[1+3g_c(1/4)
\right]+E_c(t),
\end{equation}
where $g_c(x)$ is the solution in $y$ of the equation
$$
y-(y+c)\log\left(1+\frac{y}{c}\right)= 
\log(x), 
$$
and where 
$$
E_c(t)=t\sqrt{\sum_{k=1}^d \left(\!\langle e_k,b\rangle - 
\int_{h_c(t) < \|y\|\leq 1}\langle e_k,y\rangle \nu(dy) 
+\int_{1 < \|y\|\leq h_c(t)}\langle e_k,y\rangle \nu(dy)\!\right)^2},
$$
$e_1,\ldots, e_d$ being the canonical basis of
$\bbr^d$.  
\end{thm}

\begin{rem}

(i) Note that if $X_t$ is symmetric, i.e. if $X_t\stackrel{d}{=}-X_t$, then 
$E_c(t)=0$.   

(ii) The proof of the above theorem actually shows that $3g_c(1/4)$ can be replaced by 
$g_c(1/4) + 2g_c(1/2 - t\nubar(h_c(t)))$ whenever the condition $t\nubar(h_c(t)) \leq 1/4$ is 
weakened to $t\nubar(h_c(t)) < 1/2$.

(iii) Note also that the main assumption of Theorem 1, 
namely \eqref{tnu}, is satisfied as soon as
there exists  a constant $A>0$ such that for every $R>0$, 
\begin{equation}\label{nuV}
\nubar(R)\leq A\frac{V(R)}{R^2}.  
\end{equation}
Indeed, when \eqref{nuV} holds,
choosing $c=1/4A$ ensures that $mf(X_t)$ is of order at most $h_c(t)+E_c(t)$.  
This is, in particular, 
true if $X$ is a stable vector in which case $A = (2-\alpha)/\alpha$ will do.  
In fact, in the stable case, for any $c>0$, $h_c(t) = (\sigma(S^{d-1})t/(2-\alpha)c)^{1/\alpha}$, where 
$\sigma$ is the spherical component of the corresponding stable L\'evy measure.    
In the next section, a further natural class of examples satisfying \eqref{tnu} is presented.

\end{rem}

Our next step is to study the deviations from the median. 

\begin{thm}
Under the assumptions of Theorem~1, for all $c,t>0$ such that 
$R=h_c(t)$ satisfies \eqref{nuV}, there exists $m(c,t)\in \bbr$ such 
that for all reals $x > x^\prime > 0$, the quantities 
$$\P(f(X_t)-m(c,t)\ge x) \quad {\rm and} \quad \P(f(X_t)-m(c,t)\le -x),$$ 
are upper bounded by 
\begin{equation}\label{concthm2}
Ac\frac{\nubar(x')}{\nubar(h_c(t))}
+\exp\left(\frac{x-x^\prime}{h_c(t)}-\left(\frac{x-x^\prime}{h_c(t)}+c\right)
\log\left(1+\frac{x-x^\prime}{ch_c(t)}\right)\right).
\end{equation}
In particular, if $q>0$, then for every $t>0$ such that $R=h_{q/2A}(t)$ satisfies
\eqref{nuV} and every $x>0$ such that
\begin{equation}\label{condq}
x\geq \left[1+g_{q/2A}(q/2)\right]h_{q/2A}(t), 
\end{equation}
there exists a real $m(t)$ such that
\begin{equation}\label{concq1}
\P(f(X_t)-m(t)\ge x)\leq q, 
\end{equation}
and 
\begin{equation}\label{concq2}
\P(f(X_t)-m(t)\le -x)\leq q. 
\end{equation}
\end{thm}

\begin{rem} (i) From the proof of the above theorem, it can be 
seen that $m(c,t) = \E f(Y^{(h_c(t))}_t)$, where the 
L\'evy process $Y$ is obtained by truncating the L\'evy measure of the process $X$ at  
$R=h_c(t)$, will do.  Then, taking $c= q/2A$ in $m(c,t)$ gives $m(t)$ in \eqref{concq1} and \eqref{concq2}.  
Remark also that, since $f(X_t)$ is concentrated around some value, 
this value is necessarily close to 
the median, and so $mf(X_t)$ is necessarily close to $m(c,t)$.

(ii)  In view of \eqref{eqmed}, when $X$ is symmetric, as well as 
\eqref{concq1} and \eqref{concq2} the median
and the fluctuations of $f(X_t)$ are of order $h_{1/4A}(t)$.  

(iii) It is easily seen that when $q\to 0$, $g_{q/2A}(q/2)\to 1$. So for $q$ small enough,
$\P(f(X_t)-m(t)\ge x)\leq q$ and $\P(f(X_t)-m(t)\le -x)\leq q$, 
as soon as $x\geq (2+\varepsilon)h_{q/2A}(t)$, $\varepsilon > 0$.

\end{rem}

Let us now return to the mean and let us precisely recall the result 
of Marcus and Rosi\'nski.  Let $X$ have finite 
expectation and be centered, i.e., such that  
\begin{equation}\label{zeromean}
\E(X_t)=0, 
\end{equation}
let $t>0$ and let $x_0(t)$ be the solution 
in $x$ of the equation:
\begin{equation}\label{maro}
\frac{V(x)}{x^2}+\frac{M(x)}{x}=\frac{1}{t}.  
\end{equation}
Then 
\begin{equation}\label{boundMR}
\frac{1}{4}x_0(t) \leq \E(\|X_t\|)\leq \frac{17}{8} x_0(t), 
\end{equation}
and the factor 17/8 can be replaced by 5/4 when $X$ 
is symmetric.\\

The inequality \eqref{boundMR} 
suggests that one should have fluctuations of order $x_0(t)$ 
at time $t$.  
We shall actually prove this under the following additional assumption:
There exists a constant $K$ such that for every $R>0$,
\begin{equation}\label{VM}
M(R)\leq K\frac{V(R)}{R}.  
\end{equation}
Under this last hypothesis, \eqref{maro} entails 
$$
h_{1/(1+K)}(t) \le x_0(t) \le h_1(t), 
$$
and so 
$\E\|X_t\|\asymp h_c(t)$, where
$\asymp$ means that the ratio of the two quantities is bounded above 
and below by two positive constants. 
We can now state:
\begin{thm}
Using the notation of Theorem 1, assume also that \eqref{zeromean} and \eqref{VM} hold.  
Then for all $b>0$, all $c,t>0$ such that $R=h_c(t)$
satisfies \eqref{nuV}, and for every 1-Lipschitz function $f$,
$$
\P(f(X_t)-\E f(X_t)\geq (b+cK)h_c(t))\leq
Ac+\exp\left[b - (b+c)
\log\left(1+\frac{b}{c}\right)\right].
$$
In particular, if $q>0$, then for every $t>0$, such that $R=h_{q/2A}(t)$ satisfies
\eqref{nuV} and for every $x$ such that
\begin{equation}\label{domathm3}
x\geq \left[\frac{qK}{2A}+g_{q/2A}(q/2)\right]h_{q/2A}(t), 
\end{equation}
we have
$$
\P(f(X_t)-\E f(X_t)\ge x)\leq q.   
$$

\end{thm}

\begin{rem} (i) Of course, if $X$ has finite mean but is not centered, one obtains a 
similar result by considering the L\'evy process $X_t-\E(X_t)$.

(ii) Here again, for $q$ small enough, one has $\P(f(X_t)-\E f(X_t)\ge x)\leq q$
as soon as $x\geq (1+\varepsilon)h_{q/2A}(t)$, $\varepsilon > 0$.

(iii) Above, it is clear that left tails inequalities also hold true.  For example, 
for all $x$ satisfying \eqref{domathm3}, we have: $\P(f(X_t)-m(t)\le -x)\leq q$.  

(iv)  The results on norm estimates of infinitely divisible vectors 
derived in \cite{MR} 
were used to obtain related estimates for stochastic integrals, of deterministic and, possibly, random 
predictable integrands, with respect to infinitely divisible random measures.  
Similar applications and extensions will also carry over to our settings.

\end{rem}

\section{Examples: symmetric, truncated stable processes}

In many important situations
that have been considered in the literature, the assumption \eqref{VM} 
is satisfied.  This is the case, in particular, of L\'evy processes, for which
$\nu(dx)=g(x/\|x\|)\rho(\|x\|)dx$, where $\rho$ is a function such that, say, 
$\rho(r)\asymp cr^{-\alpha-1}$, for $r$ small enough, 
while $\rho(r)\asymp cr^{-\beta-1}$, for $r$ large enough, $0< \alpha, \beta < 2$, 
or such that $\int_1^\infty r^{2}\rho(dr) <\infty$.  
Processes of this type have been introduced in physics and are also of use 
in mathematical finance, where they provide models of asset prices 
different from the usual modeling via diffusions.

Let us examine more precisely the truncated stable case.  Let $X$ be the real  
symmetric L\'evy process without Gaussian component and L\'evy measure
$$
\frac{\nu(dx)}{dx}=\frac{K}{|x|^{1+\alpha}}{\bf 1}_{\{|x|\leq M\}}, 
$$
with $K,M>0$ and $0<\alpha<2$.  Then for every $R> 0$,
$$ 
V(R)=2K\frac{\inf(R,M)^{2-\alpha}}{2-\alpha}, 
$$
and for any $c>0$, we have 
for $0\leq t\leq (2-\alpha)cM^\alpha/2K$,
$$
h_c(t)=\left(\frac{2Kt}{(2-\alpha)c}\right)^{1/\alpha}, 
$$
while for $t\geq(2-\alpha)cM^\alpha/2K$,
$$
h_c(t)= \left(\frac{2KM^{2-\alpha}t}{(2-\alpha)c}\right)^{1/2}. 
$$
Taking, say, $c= \alpha/4(2-\alpha)$, set for $0\leq t\leq \alpha M^\alpha/(8K)$,
$$
H_\alpha(t)=\left(\frac{8Kt}{\alpha}\right)^{1/\alpha}, 
$$
while for $t\geq\alpha M^\alpha/(8K)$, set 
$$
H_\alpha(t)= \left(\frac{8KM^{2-\alpha}t}{\alpha}\right)^{1/2}. 
$$
Moreover, since 
$$
\overline{\nu}(R)=\frac{2K}{\alpha}\left(\frac{1}{R^\alpha}-\frac{1}{M^\alpha}\right)
{\bf 1}_{\{R\leq M\}},   
$$
\eqref{nuV} holds with
$$
A=\frac{2-\alpha}{\alpha}.  
$$
Thus, further setting 
$$
K(\alpha)=1+3g_{\alpha/4(2-\alpha) }(1/4), 
$$
it follows from our first theorem that for every 1-Lipschitz function $f$,
$$
|mf(X_t)-f(0)|\leq K(\alpha)H_\alpha(t).  
$$
So we recovered the fact that in small time, $X$ behaves like a stable process
of index $\alpha$ while in large time, $X$ behaves like a Brownian motion. 
But furthermore, we see that the transition occurs around a time of order
$\alpha M^\alpha/K$ and we have precise bounds estimating how this transition 
happens.

Our second and third theorems also apply and give upper bounds for the fluctuations
around the median and around the mean.
For instance, choose $q>0$.
It is then easily seen that
$$
1+g_{q\alpha/2(2-\alpha)}(q/2)\leq c_\alpha, 
$$
with
$$
c_\alpha=1+\max\left(1,\frac{(1+2e)\alpha}{2(2-\alpha)}\right).  
$$
Therefore, Theorem 2 says that if $t\leq q\alpha M^\alpha/2K$, 
then there exists some $m(t)\in\bbr$ such that
$$
\P(f(X_t)-m(t)\ge x)\leq q,  
$$
as soon as
$$
x\geq
c_\alpha\left(\frac{2Kt}{q\alpha}\right)^{1/\alpha}.  
$$
On the other hand, if $t\geq q\alpha M^\alpha/2K$,
then there exists some $m^\prime(t)\in\bbr$ such that
$$
\P(f(X_t)-m^\prime(t)\ge x)\leq q,  
$$
as soon as
$$
x\geq
c_\alpha
\left(\frac{2KM^{2-\alpha}t}{q\alpha}\right)^{1/2}.  
$$
Moreover, if ones takes, $R=M$, then \eqref{nuV} is automatically satisfied. 
This amounts to taking $A=(2-\alpha)qM^\alpha/2Kt$, and so we also have
$$
\P(f(X_t)-m(t)\ge x)\leq q,  
$$
as soon as
$$
x\geq
[1+g_{Kt/(2-\alpha)M^\alpha}(q/2)]M.  
$$
To sum up, there exists some real $m(t)$ such that
if one of these two conditions holds:
\begin{itemize}
\item
$t\leq q\alpha M^\alpha/2K$ and 
$$
x\geq 
\min\left\{c_\alpha\left(\frac{2Kt}{q\alpha}\right)^{1/\alpha},
[1+g_{Kt/(2-\alpha)M^\alpha}(q/2)]M\right\}, 
$$
\item
$t\geq q\alpha M^\alpha/2K$ and 
$$
x\geq\min\left\{c_\alpha
\left(\frac{2KM^{2-\alpha}t}{q\alpha}\right)^{1/2},
[1+g_{Kt/(2-\alpha)M^\alpha}(q/2)]M\right\}, 
$$
\end{itemize}
then 
$$
\P(f(X_t)-m(t)\ge x)\leq q.  
$$
Suppose for instance that $t\leq q\alpha M^\alpha/2K$.
For $q$ not too small, the minimum is attained for the first term, and so
 the condition is $x\geq c(t/q)^{1/\alpha}$. On the other hand, for
very small $q$, the condition is $x\geq G_t(q)$ where $G_t$ can be expressed in terms
of the function $g$.

Alternatively, one can write, for $x\geq  Mc_\alpha $,
$$
\P(f(X_t)-m(t)\ge x)\leq \min \left\{\frac{C^\prime t}{x^2}, G_t(x)\right\}, 
$$
and for $x\leq Mc_\alpha $
$$
\P(f(X_t)-m(t)\ge x)\leq \min \left\{\frac{Ct}{x^\alpha}, G_t(x)\right\}, 
$$
with
$$
C=\frac{2Kc_\alpha^\alpha}{\alpha}, 
$$
$$
C^\prime=\frac{2KM^{2-\alpha}c_\alpha^2}{\alpha}, 
$$
and
$$
G_t(x)=
\exp\!\!\left[\!\!\left(\!\frac{x}{M}-1\!\right)\!-\!\left(\!\!\frac{x}{M}\!-1\!+\!\frac{Kt}{(2-\alpha)M^\alpha}\!\!\right)
\!\!\log\!\left(\!\!1\!+\!\frac{(2-\!\alpha)M^{\alpha-1}(x-\!M)}{Kt}\!\right)\!\!\right].
$$

\section{Proofs}

\subsection{Proof of Theorem 1}

Fix $t>0$, and, as in \cite{HM}, decompose $X_t$ by truncating the L\'evy measure $\nu$ at 
$R$ (to be chosen later).  Write $X_t=Y_t^{(R)}+Z_t^{(R)}$, 
where $Y^{(R)} = (Y_t^{(R)}, t\geq 0)$ and 
$Z^{(R)} = (Z_t^{(R)}, t\geq 0)$ are two independent 
L\'evy processes.  Their respective 
characteristic exponent   
${\psi_Y^{(R)}}$ and
${\psi_Z^{(R)}}$ are given, for $u\in \bbr^d$, by:   
$$
\psi_Z^{(R)}(u)=
\int_{\|y\|> R} 
\left(e^{i\langle u,y\rangle}-1\right)\nu(dy),
$$
$$\psi_Y^{(R)}(u)= i\langle u,\tilde b\rangle + \int_{\|y\|\le R} 
\left(e^{i\langle u,y\rangle}-1-
i\langle u,y\rangle{\bf 1}_{\|y\|\le 1}\right)\nu(dy),$$
with
$$
\tilde b = b-\int_{\|y\|> R} y{\bf 1}_{\|y\|\le 1}\nu(dy),
$$
where the last integral 
is understood coordinate-wise (and so is the above difference).  
Next, our global strategy is to bound $|mf(X_t)-f(0)|$ using:
\begin{eqnarray*}
|mf(X_t)-f(0)|\leq
|mf(X_t)-mf(Y_t^{(R)})|&+&|mf(Y^{(R)}_t)-\E f(Y^{(R)}_t)|\\
&+&|\E f(Y^{(R)}_t)-f(0)|.  
\end{eqnarray*}

Let us start by bounding $|mf(X_t)-mf(Y_t^{(R)})|$.  To do so, 
it is easy to check (see for instance \cite{HM}, p.1498) that 
\begin{equation}\label{boundZ}
\P(Z_t^{(R)}\ne 0)\leq t\nubar(R).  
\end{equation}
On the other hand,  
\cite{H} tells us that
$$
\P(f(Y_t^{(R)})-mf(Y_t^{(R)})\ge x)
\le H^{(R)}(x), 
$$
where
$$
H^{(R)}(x)=\exp\left(\frac{x}{2R} - \left(\frac{x}{2R}+
\frac{tV(R)}{R^2}\right)\log\left(1+\frac{Rx}{2tV(R)}\right)\right). 
$$
Next, let 
$$
I^{(R)}(y)=\sup\{x\geq 0: H^{(R)}(x)\geq y\}, 
$$
and let 
$$
P_m=\P(f(X_t)\leq mf(X_t))\geq 1/2.  
$$
Then we have (see \cite{HM} p.~1500)
\begin{equation}\label{medXY}
|mf(Y_t^{(R)})-mf(X_t)|\leq I^{(R)}(P_m-\P(Z_t^{(R)}\ne 0))\leq 
I^{(R)}(1/2-t\nubar(R)), 
\end{equation}
provided that $t\nubar(R)<1/2$.

To bound $|\E f(Y^{(R)}_t)-mf(Y^{(R)}_t)|$, we use 
the concentration inequality \cite{H}
\begin{equation}\label{expY}
\P(|f(Y_t^{(R)})-\E f(Y_t^{(R)})|\ge x^\prime)
\le 2\exp\!\!\left(\!\frac{x^\prime}{R} - \left(\!\frac{x^\prime}{R}+\frac{tV(R)}{R^2}\!\right)
\log\!\left(\!1+\frac{Rx^\prime}{tV(R)}\!\right)\!\!\right).  
\end{equation}
By the very definition of a median, and taking 
$x^\prime = |\E f(Y^{(R)}_t)-mf(Y^{(R)}_t)|$, we see that \eqref{expY} 
lead to our second estimate:  
\begin{equation}\label{expmedY}
2|\E f(Y^{(R)}_t)-mf(Y^{(R)}_t)|\leq I^{(R)}(1/4).  
\end{equation}
Finally, we bound $|\E f(Y^{(R)}_t)-f(0)|$.  
\begin{eqnarray}\label{expY'}
|\E f(Y^{(R)}_t)-f(0)|&\leq& \E\|Y^{(R)}_t\|\nonumber \\
&\leq&
\sqrt{\E\|Y^{(R)}_t\|^2}\nonumber \\
&=&
\sqrt{\!\sum_{k=1}^d\!\left(\!\!t\!\!\int_{\|y\|\leq R}\!\!\!y_k^2\nu(dy) + 
t^2\!\!\left(\!\!\tilde b_k +\!\!\int_{\|y\|\leq R}\!\!\!y_k{\bf 1}_{\|y\| > 1}\nu(dy)\!\right)^{\!\!2}\right)}\nonumber\\
&=& {\sqrt{tV(R) + \|\E Y^{(R)}_t\|^2}}.  
\end{eqnarray}  
Combining \eqref{medXY}, \eqref{expmedY} and
\eqref{expY'}, gives for any $t > 0$ and $R > 0$ such that 
$t\nubar(R)<1/2$,
$$
|mf(X_t)-f(0)|\leq I^{(R)}\!(1/2-t\nubar(h_c(t)))+ 2^{-1}I^{(R)}(1/4) + 
\sqrt{tV(\!R) + \|\E(Y^{(R)}_t\!)\|^2}.  
$$
Now, choosing $R= h_c(t)$, gives  
$$
|mf(X_t\!)-f(0)|\leq I^{(h_c(t))}(1/2-t\nubar(h_c(t)\!)\!)+ 
2^{-1}\!I^{(h_c(t))}(\!1/4) + h_c(t) + \|\E(Y^{(h_c(t))}_t\!)\|,  
$$
where $\|\E(Y^{(h_c(t))}_t)\|$ is equal to $E(t)$ given in the statement.
Finally, note that 
$$
\frac{I^{(R)}(x)}{2R} - \left(\frac{I^{(R)}(x)}{2R}
+\frac{tV(R)}{R^2} \right)\log\left(1+\frac{RI^{(R)}(x)}{2tV(R)}\right)= \log(x), 
$$ 
and so $I^{(h_c(t))}(x)=2h_c(t)g_c(x)$ with the definition of $g_c$ given in the 
statement of Theorem~1.  This concludes the proof.

\subsection{Proof of Theorem~2}

Recall the assumptions and notation as the previous subsection:
$t>0$ is fixed and $c>0$ is such that $R=h_c(t)$ satisfies \eqref{nuV}. Put
$$
m(c,t)=\E f(Y^{(R)}_t)
$$
Since $f$ is 1-Lipschitz, we have 
$|f(X_t)-f(Y^{(R)}_t)|\leq\|Z_t^{(R)}\|$.  Therefore for every
$x^\prime <x$,
\begin{equation}\label{decomp}
\P(f(X_t)-m(c,t)\ge x) \le 
\P(f(Y^{(R)}_t)-m(c,t)\ge 
x-x^\prime)+\P(\|Z_t^{(R)}\|\ge x^\prime).
\end{equation}
The first term of the above right-hand side is bounded as in  \eqref{expY}.
On the other hand, recall that $Z_t^{(R)}$ can be seen as the value
at time $t$ of a compound Poisson process $(Z_s^{(R)},s\geq 0)$. 
Therefore, if $\|Z_t^{(R)}\|\ge x^\prime$,
the process  $(Z_s^{(R)},s\geq 0)$ has at least a jump of size 
$\geq x^\prime$ before time $t$.
This implies that 
$$
\P(\|Z_t^{(R)}\|\ge x^\prime)\leq t\nubar(x^\prime).  
$$
Using \eqref{nuV}, gives 
\begin{equation}\label{boundZL}
\P(\|Z_t^{(R)}\|\ge x^\prime)\leq
Ac\frac{\nubar(x^\prime)}{\nubar(h_c(t))},
\end{equation}
which when combined with \eqref{decomp} lead to 
the inequality \eqref{concthm2} giving the first part of the theorem.  
The second part of the theorem follows by taking $x'=h_c(t)$ and $c=q/2A$.  This choice provides first 
the upper bound $q/2$ 
on \eqref{boundZL}, and moreover entails that the condition \eqref{condq} becomes  
$(x-x')/h_c(t) \geq g_c(q/2)$, leading to another upper bound $q/2$ on the rightmost term 
in \eqref{concthm2}.

\subsection{Proof of Theorem~3}

Let $c,t>0$.  Decompose $X_t$ by truncating the measure $\nu$ at 
$R=h_c(t)$.  As above, write $X_t=Y_t^{(R)}+Z_t^{(R)}$, 
where $Y_t^{(R)}$, $Z_t^{(R)}$ are two independent, 
infinitely divisible random 
vectors whose L\'evy measures are, respectively,
$t\nu(dx){\bf 1}_{\|x\|\leq R}$ and $t\nu(dx){\bf 1}_{\|x\|> R}$. Then
for every $a>cK$,
\begin{equation}\label{divide}
\P(f(X_t)-\E f(X_t)\ge ah_c(t))\le \P(f(Y_t^{(R)})-\E f(X_t)\ge ah_c(t))+
\P (Z_t^{(R)}\ne 0).
\end{equation}
Since $Z_t^{(R)}$ is a compound Poisson process, it is easy to check 
that 
\begin{equation}\label{boundZtoo}
\P(Z^{(R)}_t\ne 0)\leq t\nubar(R).  
\end{equation}
On the other hand, 
$$
\P(f(Y_t^{(R)})-\E f(X_t)\ge ah_c(t)) \leq 
\P(f(Y_t^{(R)})-\E f(Y_t^{(R)})\ge x'), 
$$ 
with
$$
x^\prime=ah_c(t)-|\E f(X_t)-\E f(Y_t^{(R)})|.
$$
To bound $x^\prime$ from below, remark that
\begin{eqnarray*}
\left|\E f(X_t)-\E f(Y_t^{(R)})\right|&=&\left|\E\left(f(Y_t^{(R)}+Z_t^{(R)})-
f(Y_t^{(R)})\right){\bf 1}_{\{Z_t^{(R)}\ne 0 \}}\right|\\
&\leq& \E\|Z_t^{(R)}\|\\
&\leq& t\int_{\|x\|>R}\|x\| \nu(dx)\\
&=& tM(R)\\
&\leq& cKR, 
\end{eqnarray*}
using both \eqref{VM} and \eqref{h} for the last inequality. Hence
$$
x^\prime\geq (a-cK)R.
$$
Moreover, \cite{H} tells us that
$$
\P(f(Y_t^{(R)})-\E f(Y_t^{(R)})\ge x')
\le \exp\left(\frac{x^\prime}{R} - \left(\frac{x^\prime}{R}+
\frac{tV(R)}{R^2}\right)\log\left(1+\frac{Rx^\prime}{tV(R)}\right)\!\!\right).  
$$
Using the fact that $R=h_c(t)$ and $x^\prime \geq (a-cK)R$, we get
$$
\P(f(Y_t^{(R)})-\E f(Y_t^{(R)})\ge x^\prime)
\le \exp\left[b - (b+c)
\log\left(1+\frac{b}{c}\right)\right], 
$$
with $b=a-cK$.
Together with \eqref{boundZtoo}, this yields the first part of Theorem~3.  
The second part follows by taking $c=q/2A$, and 
$b=g_{q/2A}(q/2)$.

\end{document}